\documentclass[12pt]{article}
\def\R{{\rm I}\!{\rm  R}}

\catcode`@=11 \@addtoreset{equation}{section} \catcode`@=12
\newtheorem{thm}{Theorem}[section]
\newtheorem{pro}[thm]{Proposition}
\newtheorem{lem}[thm]{Lemma}

\newtheorem{rema}[thm]{Remark}
\def\grad {\nabla}

\def\R{{\rm I}\!{\rm  R}}

\def\grad{\nabla} 
\def\ox{\bar x}
\def\ro {\bar r}
\def\la{\lambda}
\def\ga{\gamma}
\title{Some Liouville theorems for the $p$-Laplacian}
\author{Isabeau Birindelli \\Universit\`a di Roma ``La Sapienza''\and Fran\c coise Demengel\\University of Cergy-Pontoise}
\date{}
\begin{document}
\maketitle
\section{Introduction}
In this paper we present several Liouville type results for the $p$-Laplacian in $\R^N$. 
Let us present an example of the results here obtained.
Suppose that 
 $h$ is  a nonnegative regular function such that 
\begin{equation}\label{eqh}
h(x) = a|x|^\gamma\ {\rm for}\   |x|\ {\rm  large},\  a>0\ {\rm  and}\ \gamma>
-p. 
\end{equation}

\begin{thm}\label{th1}
1) Suppose that $N>p>1$, and 
$u\in W^{1,p}_{loc} (\R^N)\cap {\cal C} (\R^N)$ is a nonnegative weak solution 
of 
\begin{equation}\label{eq 0}
- {\rm div} (|\nabla u|^{p-2 }\nabla u) \geq h(x) u^q \;\;\mbox{in }\; \R^N
\end{equation}
with $h$ as in (\ref{eqh}).  Suppose that 
$$p-1< q\leq {(N+\ga)(p-1)\over N-p}$$
 then $u\equiv 0$.

2) Let $N\leq p$. If  $u\in W^{1,p}_{loc} (\R^N)\cap {\cal C} (\R^N)$ 
is a  weak solution 
 bounded below of 
$$-{\rm div} (|\nabla u|^{p-2 }\nabla u)\geq 0\;\;\mbox{in } \;\R^N$$
then $u$ is constant.
\end{thm}
We shall denote the $p$-Laplacian by $\Delta_p u:= {\rm div} (|\nabla u|^{p-2 }\nabla u)$.
\begin{rema}
By weak solutions of $-\Delta_p u = f$ in $\R^N$ we mean that $u\in W^{1,p}_{loc} (\R^N)$ and 
$$\int_{\R^N} |\nabla u|^{p-2} \nabla u.\nabla \varphi = \int_{\R^N} f\varphi$$
for all $\varphi\in {\cal D} (\R^N)$ (or $\varphi\in W^{1,p}(\R^N)$, compactly supported).
Similarly by weak solutions of $-\Delta_p u\geq f  (\leq f)$ we mean  that $u\in W^{1,p}_{loc} (\R^N)$ and 
$$\int_{\R^N} |\nabla u|^{p-2} \nabla u.\nabla \varphi \geq  \int_{\R^N} f\varphi (\leq  \int_{\R^N} f\varphi)$$
for all $\varphi\in {\cal D} (\R^N)$, $\varphi\geq 0$.

\end{rema}
 Non existence results for uniformly elliptic semi-linear  
 equations have been the subject of many papers. 
When $N> p=2$, the result in Theorem \ref{th1} is due to Gidas \cite{G}. 
It was extended to semi-linear equations in cones by Berestycki, Capuzzo
Dolcetta and Nirenberg in \cite {BCN1}, \cite {BCN2}. The case where the
operator is fully non linear and uniformly elliptic was treated by
Cutr\`\i\  and  Leoni  in \cite{CL}. For other non linear Liouville theorems,
see e.g.
\cite {BCDC}, \cite {SZ}....

\bigskip
We would like to remark that the first result of Theorem \ref{th1} is optimal
in the sense that for any $q>{(N+\ga)(p-1)\over N-p}$ we construct a
nonnegative solution of (\ref{eq 0}). A similar example was given in
\cite{BM} when $p=2$.

Let us also remark that the condition on $\gamma$ in (\ref{eqh}) is optimal. Indeed, for
$\gamma> -p$, Dr\'abek in \cite{Dr} has proved the existence of non trivial weak solutions in
$\R^N$ (see e.g. Theorem 4.1 of \cite{DKN}).

When treating the equation instead of the inequality, the values of $q$ for which non existence results 
hold true are not the same. Precisely  for the following equation
\begin{equation}\label{eq 01}
-\Delta_p u =   r^\gamma u^q ,\;\;u\geq 0\ \mbox{in }\; \R^N,  
\end{equation}
we prove in section 3.3 that for $p-1<q< {(N+\gamma)(p-1)+p+\gamma \over N-p}$ and $\gamma\geq 0$ 
any radial solution of (\ref{eq 01}) is identically zero. 

We believe that it should be possible to
prove that as in the case $p=2$, any solution of (\ref{eq 01}) is radial (see e.g. the work of  Damascelli and Pacella
\cite{DP},  for symmetry results for equations involving the $p$-Laplacian). Hence we
expect the non existence result to be true for any solution of (\ref{eq 01}).

When $p=2$, Gidas and Spruck \cite{GS2} have proved that for $1<q< {N+2\over N-2}$ any
solution of (\ref{eq 01}) is trivial (see also Chen and Li \cite{ChL} for a
simpler proof of this result).

Gidas and Spruck have used this to obtain a priori estimates for
solutions of the following problem:
\begin{equation}\label{eq 02}
\left\{
\begin{array}{lc}
Lu +f(x,u)=0& \mbox{in }\; \Omega\\
u=\phi & \mbox{on }\; \partial\Omega
\end{array}
\right.
\end{equation}
where $L$ is a second order uniformly elliptic operator and $f$ satisfies some
growth conditions. This is done through a blow up argument (see also
\cite{BCN2}).

Analogously Theorem \ref{th1} constitutes the first step to obtain a priori
estimates for reaction diffusion equations involving $p$-Laplacian type
operators in bounded domains, this will be the object of future works.

The other Liouville type theorems here enclosed concern the following equation

\begin{equation}\label{eq 03}
\Delta_p u \geq h(x) u^q \;\;\mbox{in }\; \R^N
\end{equation}
for $h$ as in (\ref{eqh}).
We prove that if 
$p-1< q$ and $\gamma>-p$
then $u\equiv 0$.

For completeness sake, we begin this paper by proving that bounded $p$-harmonic
functions in $\R^N$ are constant. This result is probably known, but since we
have not found it in the litterature, we include its proof here.

\section{$p$-harmonic functions}

Recall that $u$ is $p$-harmonic in a domain $\Omega$ if  $u\in W^{1,p}_{loc} (\Omega)$ and 
\begin{equation}\label{eqph}
-\Delta_p u = -{\rm div} (|\nabla u|^{p-2} \nabla u )
= 0\ {\rm in} \ \Omega .
\end{equation}
\begin{thm}\label{thph}
Suppose that $u$ is a bounded $p$-harmonic function in $\R^N$, then $u$ is constant.
\end{thm}
Before giving the proof of Theorem \ref{thph} we
we need three lemma:

\begin{lem}\label{lemma 1}
Suppose that $u$ is $p$-harmonic and belongs to
$W^{1,p}_{loc} ({\bf R}^N)$. Then, there exists some
constant $C$ such that for all
$R$ and 
$\sigma$ $>0$
$$\int_{B(0,R)} |\nabla u|^p\leq c{1\over
\sigma^p} \int_{B(0,R+\sigma)} |u|^p.$$
\end{lem}

\begin{lem}\label{Lemma 2}
Suppose that $u$ is and
$p$-harmonic then, there exists some positive
constant $C$ such that for all $\lambda\geq 1$, $R$ and
$\sigma$
$$ \int_{B(0,R)}
|\nabla(|\nabla u|^{p+\lambda\over 2})|^2\leq
{C\over \sigma^2} \int_{B(0,R+\sigma)}|\nabla
u|^{p+\lambda}.$$
In particular for all $k$ and $n\in {\bf N}$
$$\int_{B(0,R)} |\nabla (|\nabla u|^{pk^{n-1}\over
2})|^2\leq {C\over \sigma^2}
\int_{B(0,R+\sigma)} |\nabla u|^{pk^{n-1}}.$$
 \end{lem}

\begin{lem}\label{Lemma 3}
Suppose that $\phi_n$ is some sequence of
positive numbers satisfying for some constants
$c>0$ and $k>0$
$$\phi_n\leq c^n\phi_{n-1}^k$$
Then 
$$\phi_n^{k^{-n}}\leq c^{k\over (k-1)^2} \phi_0.$$
\end{lem}
{\bf Proof of Lemma \ref{Lemma 3}}

Let us start by noticing that 
\begin{eqnarray*}
\phi_n&\leq &c^n\phi_{n-1}^k \\
&\leq& c^{n+(k)(n-1)} \phi_{n-2}^{k^2}\\
&\leq& c^{\sum_{p=0}^{n-1} (n-p)k^p}\phi_0^{k^n}
\end{eqnarray*}
One can easily compute
$$\sum_{p=0}^{n-1} (n-p)k^p= (n+1){k^n-1\over k-1}-
{d\over dk}\left({k^n-1\over
(k-1)^2}\right)$$
and then
$$\phi_n^{k^{-n}}\leq c^{k^{n+1}\over
(k-1)k^n}\phi_0
\leq
c^{k\over (k-1)^2}\phi_0.$$
{\bf Proof of Lemma \ref{lemma 1}}

First let us remark that by the regularity results of Tolksdorff \cite{To} $u$ is ${\cal
 C}^1(\R^N)$ and $u\in W^{2,p}_{loc}$ for $p<2$, $u\in W^{2,2}_{loc}$ for $p\geq 2$. Multiply
equation (\ref{eqph}) by
$u\zeta^p$, where
$\zeta$ is some smooth function which equals one
on $B(0,R)$ and is zero outside $B(0,R+\sigma)$
$0\leq \zeta\leq 1$. One has then some
universal constant such that $|\nabla
\zeta|_\infty\leq {c\over
\sigma}$. One obtains
$$\int |\nabla u|^p\zeta^p= -p\int |\nabla
u|^{p-2} \nabla u.(\nabla \zeta)\zeta^{p-1}u
\leq |(\nabla u) \zeta|_p^{p-1} (\int u^p|\nabla
\zeta|^p)^{1\over p}.$$
Finally
$$\int_{B(0,R)} |\nabla u|^p\leq \int |(\nabla u)
\zeta|^p\leq {c\over \sigma^p}
\int_{B(0,R+\sigma)-B(0,R)}|u|^p.$$
In particular, if $u$ is bounded
$$\int_{B(0,R)}|\nabla u|^p\leq
C{(R+\sigma)^N-R^N\over
\sigma^p}.$$
For $R=\sigma$ one obtains
$$\int_{B(0,R)}|\nabla u|^p\leq {2^NC^N} R^{N-p}.$$

\noindent {\bf Proof of Lemma \ref{Lemma 2}}

We use  formal calculations which can be justified, as it is done in
\cite{To}.
 Differentiate equation (\ref{eqph}) with respect to
$x_k$ and multiply  by  $|u_{,k}|^\lambda
u_{,k}\zeta^2$, integrating we obtain
\begin{eqnarray*}
\int(|\nabla u|^{p-2}u_{,j})_{,k}
(|u_{,k}|^\lambda u_{,k}\zeta^2)_{,j}&= &
(\lambda+1)\int|\nabla u|^{p-2}
|u_{,k}|^\lambda(u_{,jk}^2)\zeta^2
+\\
&+&(p-2)(\lambda+1)\int|\nabla u|^{p-4} (\nabla
u.\nabla u_{,k})^2\zeta^2+\\
&+&\int|\nabla u|^{p-2}
|u_{,k}|^\lambda u_{,k}2\zeta\zeta_{,j} u_{,kj}\\
&+&(p-2)\int |\nabla u|^{p-4}
(\nabla u.\nabla u_{,k}) 2\zeta \zeta_{,j}| u{,k}|^\lambda u_{,k}. 
\end{eqnarray*}
Let us observe that
$$|\int|\nabla u|^{p-4} (\nabla
u.\nabla u_{,k})^2\zeta^2|\leq \int|\nabla u|^{p-2}
|u_{,k}|^\lambda(u_{,jk}^2)\zeta^2$$
by Schwartz inequality. 
Using $|p-2|\leq 1$ if $p<2$ and if not 

\noindent $(p-2)(\lambda+1)(\int|\nabla u|^{p-4} (\nabla
u.\nabla u_{,k})^2\zeta^2)\geq 0$, summing
over $k$ and  using once more Schwartz inequality
one obtains
\begin{eqnarray*}
(\lambda+1)(p-1)\int|\nabla
u|^{p-2} |\nabla u|^\lambda|\nabla \nabla
u|^2\zeta^2 &\leq&2(p-1) \int|\nabla
u|^{p+\lambda-1}|\nabla \nabla u
|\zeta|\nabla \zeta|\\ &\leq& {(\lambda+1)(p-1)\over
2}\int |\nabla u|^{p+\lambda-2}|\nabla
\nabla u|^2+\\
&+&{C\over
\lambda+1}\int |\nabla u|^{p+\lambda}|\nabla
\zeta|^2.
\end{eqnarray*}
Finally
$$\int_{B(0,R)}|\nabla
u|^{p-2+\lambda} |\nabla \nabla u|^2\leq
{C\over (\lambda+1)^2} {1\over \sigma^2}
\int_{B(0,R+\sigma)} |\nabla u|^{p+\lambda},$$
i.e.
$${4\over (p+\lambda)^2}\int|\nabla (|\nabla
u|^{p+\lambda\over 2})|^2\leq {C\over
(\lambda+1)^2 \sigma^2
}\int_{B(0,R+\sigma)}|\nabla u|^p$$
which implies , iterating and beginning with
$\lambda=0$,
$$\int_{B(0,R)}|\nabla (|\nabla
u|^{pk^{n-1}\over 2})|^2\leq {C\over
 \sigma^2
}\int_{B(0,R+\sigma)}|\nabla u|^{pk^{n-1}}.$$
This concludes the proof of Lemma \ref{Lemma 3}

\noindent{\bf Proof of Theorem \ref{thph}.}
Let us recall that according to some
Poincar\'e's inequality
and Sobolev embedding one has that for all $k\leq
{N\over N-p}$, there exists some constant $C_N$
which depends only on $N$ and $p$ such that for
all $R$ and $w\in H^1_{loc}({\bf R}^N)$
$$\left({1\over
R^N}\int_{B(0,R)}|w|^{2k}\right)^{1\over k}\leq
C_N\left({1\over R^{N-2}}\int_{B(0,R)}|\nabla
w|^2+{1\over R^N} \int_{B(0,R)}|w|^2\right).$$
Let us choose   
$w= |\nabla u|^{pk^n\over 2}$ in this inequality in order to  obtain
$$\left({1\over
R^N}\int_{B(0,R)}|\nabla u|^{pk^n}\right)^{1\over
k}\leq C_N\left({1\over
R^{N-2}}\int_{B(0,R)}|\nabla (|\nabla
u|^{pk^n\over 2})|^2+{1\over R^N}
\int_{B(0,R)}|\nabla u|^{pk^{n-1}}\right).$$
We are in a position to apply  Lemma \ref{Lemma 2}, and the previous inequality becomes
\begin{eqnarray}\label{eqx}
\left({1\over
R^N}\int_{B(0,R)}|\nabla u|^{pk^n}\right)^{1\over
k}&\leq& C^\prime_N ({1\over R^N})(1+{R^2\over
\sigma^2})\int_{B(0,R+\sigma)} |\nabla
u|^{pk^{n-1}} \\
&\leq& C^\prime_N (1+{R^2\over
\sigma^2})({R+\sigma\over R})^N {1\over
(R+\sigma)^N}\int_{B(0,R+\sigma)} |\nabla
u|^{pk^{n-1}}.\nonumber
\end{eqnarray}
We define $\phi_n
=\left({1\over
\rho_n^N}\int_{B(0,\rho_n)}|\nabla
u|^{pk^n}\right)^{1\over k}$
with $\rho_n= r(1+2^{-n})$, which satisfies
$\rho_{n-1}\leq 2\rho_n$. 
 With $R= \rho_n$ and $\sigma=
\rho_{n-1}-\rho_n$, (\ref{eqx}) becomes 
$$\phi_n\leq C(1+\left({\rho_n\over
\rho_{n-1}-\rho_n}\right)^2) \phi_{n-1}^k$$
Using 
$$1+\left({\rho_n\over
\rho_{n-1}-\rho_n}\right)^2= 1+(2^n+1)^2\leq 5^n,$$
one finally has 
$$\phi_n\leq C5^n \phi_{n-1}^k.$$
Using Lemma \ref{lemma 1}, one gets that 
$$\overline{\rm lim}_{n\rightarrow\infty}
\phi_n^{1\over k^n} \leq C\phi_0= C{1\over
(2r)^N}\int_{B(0,2r)} |\nabla u|^p.$$
The left hand side tends to $\sup_{B(0,r)}
|\nabla u|^p$ and the right hand side is less
than a constant multiplied by $(2r)^{N-p-N}$ by
using Lemma \ref{lemma 1}. Taking $r$ large enough one gets
that $\nabla u=0$, which ends the proof.

\section {The positive semi-linear case}

\subsection{The inequation}
When $N>p$ our main non-existence result in this section is
the following

\begin{thm}\label{lnl}
Suppose that $N>p>1$. Let  $u\in W^{1,p}_{loc} (\R^N)\cap {\cal C} (\R^N)$ be a nonnegative weak
solution  of 
\begin{equation}\label{eq}
-\Delta_p u \geq h(x) u^q \;\;\mbox{in }\; \R^N,
\end{equation}
with $h$ satisfying (\ref{eqh}).
 Suppose that 
$p-1< q\leq {(N+\ga)(p-1)\over N-p}$,
then $u\equiv 0$.
\end{thm}
The proof is inspired by the one given in \cite{CL}, where the authors treat fully nonlinear
strictly elliptic equations. 

 Let us start by one remark and two propositions.
\begin{rema} \label{r1} The following comparison result holds true:
let $u$ and $\phi$ satisfy $u, \phi\in W^{1,p}(\Omega)$
\[
\left \{\begin{array}{lc}
-\Delta_pu\geq -\Delta_p\phi =0 & \mbox{in} \; \Omega\\
u\geq \phi   & \mbox{on} \; \partial\Omega
\end{array}
\right.
\]
then $u\geq \phi$ in $\Omega$.
This is a standard result and it is easy to see for example by 
multiplying  $-\Delta_p u+\Delta_p \phi$
by $(\phi-u)^+$. 
\end{rema}
\begin{pro}\label{prop10}
Let $\Omega$ be an open set in $\R^N$, and let $f\in {\cal C}(\overline{\Omega})$. Suppose that
$u\in W^{1,p}_{loc} (\Omega)\cap {\cal C}(\overline{\Omega}) $ is a weak solution
of
$-\Delta_p u\geq f$ in $\Omega$. Then,
  if $x_0\in \Omega$,and  $\varphi\in {\cal C}^2(\Omega)\cap {\cal C} (\overline{\Omega})$,are
such that 
 $$\nabla \varphi(x_0)\neq 0,\  u(x_0)- \varphi(x_0) = \inf_{y\in \Omega} u(y)-\varphi(y)$$ 
 then 
$$-\Delta_p \varphi(x_0)\geq f(x_0)$$
\end{pro}
This proof is inspired by Juutinen \cite{Juu}. 

\noindent{\bf Proof.} Without loss of generality we can suppose that $u(x_0)=\varphi(x_0)$.

Let us note first that it is sufficient to prove that the property holds for every $\varphi$ such
that 
$\varphi(y)< u(y)$ for all $y\neq x_0$ in a sufficiently small neighbourhood of $x_0$. Indeed,
suppose that the property holds for such functions then taking $\varphi_\epsilon (y) =
\varphi(y)-\epsilon {|y-x_0|^4}$ and letting
$\epsilon$ go to zero, one obtains the result for every $\varphi$. 

Suppose by contradiction that there exists some
$x_0\in \Omega$ and some ${\cal C}^2 $ function $\varphi$ such that $\nabla \varphi(x_0)\neq 0$, 
$\varphi(x_0)= u(x_0)$ and $\varphi(y)< u(y)$ on some ball $B(x_0,r)$ and $-\Delta_p
\varphi(x_0)<f(x_0)$. By continuity, one can choose  
$r$ sufficiently small such that 
$\nabla \varphi(y)\neq 0$ , as well as
$$-\Delta_p \varphi(y)<f(y),$$
for all $y\in B(x_0,r)$.
 Let $m = \inf_{|x-x_0|=r} (u(x)-\varphi(x))>0$, and define 
$$\bar\varphi = \varphi+{m\over 2}.$$
One has $-\Delta_p \bar\varphi<f$ in $B(x_0,r)$ and $\bar\varphi\leq u$ on $\partial B(x_0,r)$.

Using the comparison principle one gets that $\bar\varphi\leq u$ in the ball which contradicts
$\bar\varphi(x_0) = \varphi(x_0)+{m\over 2}> u(x_0)$.
This ends the proof of Proposition \ref{prop10}.

\bigskip

Finally let us recall that if $v$ is radial i.e. $v(x)=V(|x|)\equiv V(r)$ for some function $V$
 and $V$ is ${\cal C}^2$,  then if $x$ is such that $V^\prime (|x|)\neq 0$
$$\Delta_p v(x) =|V^\prime(r)|^{p-2}\left((p-1)V^{\prime\prime}(r)+{N-1\over
r}V^\prime(r)\right).$$ 
Hence  for any constants $C_1$ and $C_2$ if $N\neq p$
and for $\lambda={p-N\over p-1}$ the function $\phi(x)=C_2|x|^{\lambda}+ C_1$
satisfies $\Delta_p \phi=0$ for $x\neq 0$.

Before giving the proof of Theorem \ref{lnl} let us define 
$m(r)=\inf_{x\in B_r}u(x)$ and prove the following Hadamard type inequality 
\begin{pro} \label{p1} Let $N\neq p$. Suppose that $-\Delta_p u\geq 0$ and $u\geq 0$. Let
$\la = {p-N\over p-1}$.  For any $0<r_1<r<r_2$ :
\begin{equation}\label{i2}
m(r)\geq {m(r_1)(r^\la-r_2^\la)+m(r_2)(r_1^\la-r^\la)\over r_1^\la-r_2^\la}.
\end{equation}
Let $N=p$ then
\begin{equation}\label{i2b}
m(r)\geq {m(r_1)\log({r\over r_2})+m(r_2)\log({r_1\over r})\over\log({ r_1\over	r_2})}.
\end{equation}
\end{pro}

\noindent{\bf Proof:} Let $N\neq p$. Let  $0<r_1<r_2$.
Let us consider $\phi(r)=C_2r^{\lambda }+ C_1$ with $C_2$ and $C_1 $ such
that
$\phi(r_1)=m(r_1)$ and $\phi(r_2)=m(r_2)$. It is easy
 to see that
$$\phi(r)={m(r_2)(r^\la-r_1^\la)+m(r_1)(r_2^\la-r^\la)\over r_2^\la-r_1^\la}.$$
Obviously $\phi>0$ and
for $i=1$ and $i=2$, $u(x)\geq m(r_i)=\phi(r_i)$ for $x\in \partial B_{r_i}$, 
hence $u$ and $\phi$ satisfy the conditions of Remark \ref{r1}.
and $u(x)\geq \phi(|x|)$ in $B_{r_2}\setminus B_{r_1}$. Taking the infimum
we obtain that $\inf_{|x|=r} u(x)\geq \phi(r)$ for $r\in [r_1,r_2]$. By the 
minimum principle $m(r)=\inf_{|x|=r} u(x)$. This ends the proof of the first part of
proposition \ref{p1}.

For $N=p$ consider 
$$\psi(r)= {m(r_1)\log({r\over r_2})+m(r_2)\log({r_1\over r})\over\log({ r_1\over r_2})}.
$$
Remark that $\Delta_N\psi=0$ and $\psi(r_1)=m(r_1)$ and $\psi(r_2)=m(r_2)$. Now
proceed as above.

\bigskip

\begin{rema}\label{r3}
Clearly if $\lambda<0$ i.e. $p< N$, then $g(r):=m(r)r^{-\la}$ is an 
increasing function. Just observe that $r_1^\la-r^\la\geq 0$ and  let 
$r_2$ tend to $+ \infty$  in (\ref{i2}) and one obtains for $r\geq r_1$:
$$ m(r)\geq { m(r_1)r^\la\over r_1^\la}.$$

\end{rema}

\noindent {\bf Proof of Theorem \ref{lnl}.} 
We suppose by contradiction that $u\not\equiv 0$ in $\R^n$,
 but since $u\geq 0$ by the strict maximum principle of Vasquez \cite{V} we get that
$u>0$.

Let $0<r_1<R$, define 
$g(r)=m(r_1)\left\{1-{[(r-r_1)^+]^{k+1}\over (R-r_1)^{k+1}}\right\}$ with $k$ such that
$$k\geq 3 \;\;\mbox{and}\;\; {1\over k}< p-1.$$  Let $\zeta(x)=g(|x|)$. Clearly for $|x|<r_1$, $u(x)> m(r_1) =\zeta(x)$
 while for 
$|x|\geq R$, $\zeta(x)\leq 0 <u(x)$. 
On the other hand there exists $\tilde x$ such that $|\tilde x|=r_1$ and 
$u(\tilde x)=\zeta(\tilde x)$. Hence the minimum of $u(x)-\zeta(x)$ occurs for 
some $\ox$ such that $|\ox|=\bar r$ with $r_1\leq \bar r<R$.

Let $|x|=r$, it is an easy computation to see that for
$r\geq r_1$ 
\begin{eqnarray}\label{i3}
& &-\Delta_p\zeta(x)= \\
&=&\left({(k+1)m(r_1)\over (R-r_1)^{k+1}}\right)^{(p-1)}
\left[2(p-1)+(N-1){(r-r_1)^+\over r}\right]((r-r_1)^+)^{kp-(k+1)}.\nonumber
\end{eqnarray}
Clearly with our choice of $k$, $kp-(k+1)>0$ and hence, for $|x|=r_1$, $-\Delta_p\zeta (x)=0$
 while, of course, $\grad \zeta (x)=0$.

Now we have two cases. 

First case $\bar r=r_1$. This implies
$$u(\bar x)-m(r_1)= u(\bar x)-\zeta(\bar x)\leq u(x)-\zeta(x)$$
for all $x$. In particular choosing $  x=\tilde x$, one gets 
$$u(\bar x)-m(r_1)\leq u(\tilde x)-\zeta(\tilde x)=0.$$
Finally
$$u(\bar x)= m(r_1)$$
and $\bar x$ is a mimimum for $u$ on $B(0,r_1)$. Since $-\Delta_p u\geq 0$,
Hopf's principle as stated in Vasquez \cite{V} implies that
$\nabla u(\bar x)\neq 0$. On the other hand $\nabla u (\bar
x)= \nabla \zeta (\bar x)= 0$, a contradiction.

Second case : $r_1< \bar r< R$. Now $\nabla \zeta (\bar x)\neq 0$, and using
Proposition \ref{prop10} one has 
$$h(\bar x)u^q(\bar x)\leq -\Delta_p \zeta (\bar x).$$
We choose $r_1$ and $R$ sufficiently large in order that $h(x) = a|x|^\gamma$ for $|x|\geq \min
(r_1, {R\over 2})$.  Combining this with (\ref{i3}), we obtain

$$a\ro^\ga m(\ro)^q\leq a\ro^\ga u^q(\ox)\leq
(k+1)^{(p-1)}(N+2p-3)m(r_1)^{(p-1)}(R-r_1)^{-p}.$$ Since $m$ is decreasing we have
obtained for some constant $C>0$
$$ m(R)\leq C m(r_1)^{(p-1)\over q}\ro^{-\ga\over q}(R-r_1)^{-p\over q}.$$

Now we choose $r_1={R\over 2}$, we use Remark \ref{r3} and  the previous 
inequality becomes
$$ m(R)\leq C m(R)^{(p-1)\over q}R^{-p-\ga\over q},$$
which implies 
\begin{equation}\label{i5}
m(R)R^{-\la }\leq CR^{-\la-{p+\ga\over q-p+1}}.
\end{equation}
Clearly $-\la-{p+\ga\over q-p+1}={N-p\over p-1}-{p+\ga \over q-p+1}\leq 0$ when 
$ q\leq {(N+\ga) (p-1)\over N-p}$. 

If $q< {(N+\ga) (p-1)\over N-p}$ we have reached 
a contradiction since the right hand side of (\ref{i5})
tends  to zero for 
$R\rightarrow +\infty$ while the left
 hand side is an increasing positive function as seen in
Remark \ref{r3}.

\bigskip

\bigskip

We now treat the case $q={(N+\ga)(p-1)\over N-p}$. Let us remark that for this choice 
of $q$ we have that for some $C_1>0$, $c>0$ and $r>r_1>0$, with $r_1$ large enough :
 \begin{equation}\label{if1}
-\Delta_p u\geq ar^\gamma u^q\geq C_1r^{-N}\;\mbox{ since }
 \; m(r)\leq c r^{p-n\over p-1}.
\end{equation}

 We choose
$\psi(x) = g(|x|)$ with
$$g(r) = \gamma_1 r^{p-N\over p-1} \log^\beta r+\gamma_2$$
where $\gamma_1$ and $\gamma_2$ are two positive constants such that for some $r_1>1$ and some
$r_2>r_1$: 
$$m(r_2) = g(r_2),$$
$$m(r_1)\geq g(r_1),$$
while $\beta$ is a positive constant to be chosen later.
It is easy to see that 
\begin{eqnarray*}
\Delta_p \psi&= &|\gamma_1|^{p-1} r^{-N} \left\vert{p-N\over p-1}\log^\beta
r +\beta \log^{\beta-1} r\right\vert^{p-2}\cdot \\
 &\cdot & \left[(p-1)\beta(\beta-1)\log^{\beta-2} r-\beta(3N-2p-2)\log^{\beta-1} r\right]
\end{eqnarray*}
Suppose now that $p>2$, and choose $0<\beta<{1\over p-1}<1$, then there exists $C>0$ such that
$$
\Delta_p \psi\geq -|\gamma_1|^{p-1}Cr^{-N}(\log r)^{\beta(p-1)-1}\geq -|\gamma_1|^{p-1}Cr^{-N}(\log r_1)^{\beta(p-1)-1}.
$$

On the other hand for $p\leq 2$ we can choose $\beta=1$ and a calculation similar to the one above implies that 
$$\Delta_p \psi\geq -c|\gamma_1|^{p-1}r^{-N} \left(\log r_1\right)^{p-2}.$$
In both cases we can  choose $\gamma_1$ small enough to get  
$$\Delta_p \psi \geq -C_1r^{-N}\geq \Delta_p u.$$
Since $u\geq \psi$ on the boundary of $B_{r_2}\setminus B_{r_1}$, one obtains by
the comparison  principle (Remark \ref{r1}) that $u\geq \psi$ everywhere in $B_{r_2}\setminus
B_{r_1}$. When
$r_2$ goes to infinity it is easy to see that $\gamma_2\rightarrow 0$, and we obtain
$$u(x)\geq c |x|^{p-N\over p-1} \log^\beta |x|,$$
for $|x|\geq r_1$.
This implies that 
$$m(r)\geq cr^{p-N\over p-1} \log r$$
for $r> r_1$. We have reach a contradiction since   
$$m(r)\leq Cr^{p-N\over p-1}.$$
Hence $u\equiv 0$. This  concludes the proof of Theorem \ref{lnl}.

\bigskip

We treat now the case $N\leq p$ where the result is much stronger. 

\begin{thm}\label{Np}
Let $N\leq p$. If $u\in W^{1,p}_{loc}(\R^N)\cap {\cal C} (\R^N)$ is   bounded below and is a weak
solution of 
$$-\Delta_p u\geq 0\;\;\mbox{in } \;\R^N$$
then $u$ is constant.
\end{thm}

\begin{rema}
For $N\leq p$, for any $q>0$ and for any  nonnegative $h$, if $u\in W^{1,p}_{loc}(\R^N)\cap {\cal
C} (\R^N)$ is a weak solution  of
$$-\Delta_p u\geq h(x)  u^q \;\;\mbox{in } \;\R^N$$
then $u\equiv 0$. 
\end{rema}
{\bf Proof of Theorem \ref{Np}.} Without loss of generality we can suppose that $u\geq 0$. First
we will consider
$N<p$. Let
$m(r)=\inf_{x\in B_r(0)}u(x)$. From Proposition \ref{p1} we know that for $0<r_1<r<r_2$
\begin{equation}\label{i7}
m(r)\geq
{m(r_1)(r_2^\la-r^\la)+m(r_2)(r^\la-r_1^{\la})\over
r_2^\la-r_1^\la,}
\end{equation}
where $\la={p-N\over p-1}>0$.

 If we let $r_2\rightarrow
+\infty$ inequality (\ref{i7}) becomes
\begin{equation}\label{i8}
m(r)\geq m(r_1).
\end{equation}
But of course $m(r)$ is decreasing hence (\ref{i8}) implies that
$m(r)$ is constant i.e. $m(r)=m(0)=u(0)$ for any $r>0$.
Clearly this can be repeated with balls centered in any point of $\R^N$. Hence $u$ is constant.

For the case $N=p$ just use inequality (\ref{i2b}) in Proposition \ref{p1} and proceed as above.

This concludes the proof of Theorem \ref{Np}.

\subsection{Counterexample}
We are going to show that for $N>p$, for $\gamma\geq 0$ and for  any $q>{(N+\gamma)(p-1)\over
N-p}$  there exists a non-negative function $u$ such that 
$$-\Delta_p u\geq r^\gamma u^q \;\;\mbox{in }\;\R^N$$
hence proving that ${(N+\gamma)(p-1)\over N-p}$ is an optimal upper bound  for $q$ in Theorem \ref{lnl}.

Indeed consider $g(r)=C(1+r)^{-\alpha}$ with $\alpha$ and $C$ two positive 
constants to be determined.
Clearly $\Gamma(x)=g(|x|)$ satisfies
\begin{eqnarray*}
-\Delta_p\Gamma &=&C^{p-1}\alpha^{p-1}(1+r)^{-(\alpha+1)(p-2)}
[-(\alpha+1)(p-1)(1+r)^{-(\alpha+2)}+\\
&+&{(N-1)\over r}
(1+r)^{-(\alpha+1)}]\\
&\geq& C^{p-1}\alpha^{p-1}(1+r)^{-\alpha(p-1)-p}[N-1-(\alpha+1)(p-1)]
\end{eqnarray*}
with $r=|x|$.

Now let $\epsilon>0$ such that $q={(N+\gamma-\epsilon)(p-1)\over (N-p-\epsilon)}$
and let $\alpha={N-p-\epsilon \over p-1}$. Clearly we have 
$\alpha (p-1)+p+\gamma=N+\gamma-\epsilon=\alpha q$.  
Furthermore $N-1-(\alpha+1)(p-1)=N-p -\alpha(p-1)=\epsilon+\gamma>0$.
Hence choosing $C$ such that $C^{p-1}\alpha^{p-1}(\epsilon+\gamma)=C^q$ we obtain that
$\Gamma(x)$ satisfies
$$-\Delta_p \Gamma\geq C^q (1+r)^\gamma (1+r)^{-\alpha
(p-1)-p-\gamma}\geq r^\gamma\Gamma^q \;\;\mbox{in }\;\R^N.$$

\subsection{The equation}
In this section we are interested in studying non-existence results concerning the equation. Clearly in view of Theorem \ref{Np}, we are only interested in the case $N>p$.:
\begin{thm}
Suppose that $u\in W^{1,p}_{\rm loc}(\R^N)$ is nonnegative and satisfies
\begin{equation}\label{st1}
-\Delta_p u = r^\gamma u^q,
\end{equation}
for some $\gamma \geq 0$. 
If $$p-1<q\leq {(N+\gamma)(p-1)+p+\gamma\over N-p}$$
and  $u$ is radial then  $u\equiv 0$.
\end{thm}

\begin{rema}
One can get the same result for $-\Delta_p u = Cr^\gamma u^q$ by considering $u$
 multiplied by some convenient constant.
\end{rema}

The proof given here is similar to the one given by Caffarelli, Gidas and Spruck in \cite{CGS}.

\noindent {\bf Proof.}

It is sufficient to consider the case 
$q\geq {(N+\gamma)(p-1)\over N-p}$, since the other cases are
proved in Theorem \ref{lnl}.

If $u$ is a radial solution and satisfies (\ref{st1}) in a weak sense, then it is not difficult to
see that it satisfies in the weak sense 
$$-(r^{N-1} |u^\prime|^{p-2} u^\prime)^\prime = r^{N-1+\gamma} u^q.$$
Integrating  between $0$ and $r$, one has 
$$r^{N-1} |u^\prime|^{p-2} u^\prime = -\int_0^r s^{N-1+\gamma} u^q (s) ds.$$
Since $u^\prime<0$, $u$ is decreasing and then, 
$$r^{N-1} |u^\prime|^{p-2} u^\prime \leq -u(r)^q {r^{N+\gamma}\over N+\gamma}.$$
Hence
$$u^\prime u^{-q\over p-1} \leq -cr^{1+\gamma\over p-1}$$
and integrating one gets 
$$u(r)\leq C r^{\gamma+p\over p-1-q}.$$
Coming back to the equation one obtains 
$$r^{N-1} |u^\prime|^{p-1} = \int_0^r s^{N-1+\gamma} u^q (s) ds \leq C\int_0^r
s^{N-1+\gamma } s^{(\gamma+p)q\over (p-1-q)}ds.$$
Clearly
$N+\gamma+{(\gamma+p)q\over p-1-q}\geq 0$ when 
$q\geq {(N+\gamma)(p-1)\over N-p}$ and therefore 
$$|u^\prime(r)|^{p-1}\leq C r^{\gamma+  {(\gamma+p)q\over p-1-q}+1}$$
and then
$$|u^\prime |\leq Cr^{(\gamma+q+1)\over p-1-q}.$$
In order to conclude, we need to use Pohozaiev identity:
$$(N-p)\int_B |\nabla u|^p+ p\int_{\partial B} \sigma.n (\nabla u.x)=
p\int_B\Delta_p u (\nabla u.x)+\int_{\partial B} |\nabla u|^p (x.\vec
n)$$
here $\sigma = |\nabla u|^{p-2} \nabla u$ and $B=B(0,R)$. 
>From the equation we know that 
$$\int_B |\nabla u|^p-\int_{\partial B} (\sigma.\vec n) u = \int_B r^\gamma
u^{q+1}$$
and then
\begin{eqnarray}
(N-p)\left(\int_B r^\gamma u^{q+1}\right.&+&\left.\int_{\partial B} \sigma.\vec n u\right) + 
p\int_{\partial B} (\sigma.\vec n ) (\nabla u.x)\nonumber\\
& = & -p\int_B r^\gamma u^q
(\nabla u.x)+\int_{\partial B} |\nabla u|^p x.\vec n.
\end{eqnarray}
Using the fact that $u$ is radial, for $\omega_n=|B_1|$ one gets
\begin{eqnarray*}
{1\over \omega_n}\int_{B_R} r^\gamma u^q \nabla u.x dx&=& \int_0^R r^{\gamma+N} u^q(r)
u^\prime(r) dr \\
&=& 
\int_0^R r^{\gamma+N} ({u^{q+1}(r)\over q+1})^\prime dr\\
& =& -{\gamma+N\over q+1}
\int_0^R r^{\gamma+N-1} u^{q+1}+{1\over q+1}R^{\gamma+N} u^{q+1}(R).
\end{eqnarray*}
We have
finally obtained
\begin{eqnarray*}
\left(N-p-{(\gamma+N)p\over q+1}\right) \int_0^R r^{\gamma+N-1} u^{q+1}
dr&= &(N-p) |u^\prime|^{p-1}u r^{N-1} + (1-p) |u^\prime|^p r^N \\
&-&{p\over
q+1}r^{\gamma+N} u^{q+1}.
\end{eqnarray*}
Let us note that since $q< {(N+\gamma)(p)+p-N\over N-p}$,
one has 
$${(\gamma+N)p\over q+1}+p-N>0.$$ Moreover the estimates on $u$ and $u^\prime$
imply that the terms $|u^\prime|^{p-1}u(R) R^{N-1}$, $ |u^\prime|^p(R)
R^N$ and $R^{\gamma+N} u^{q+1}(R)$ behave respectively as 
$R^{N-1+{\gamma+p\over p-1-q}+{(\gamma+q+1)(p-1)\over
p-1-q}}$, 
$R^{\gamma+N+{\gamma+p\over p-1-q}(q+1)}$ and $R^{N-p\left({\gamma+q+1\over
q-p+1}\right)}$. All the exponents are negative, and then 
$\int_0^R r^{\gamma+N-1} u^{q+1}
dr\rightarrow 0$ when $R\rightarrow +\infty$, hence $u\equiv 0$.  
This conclude the proof.

\section{The negative semi-linear case}
In this section we prove analogous results for inequations in which $\Delta_p$ is replaced by
$-\Delta_p$, in this case no upper bound for $q$ is required.  
\begin{thm}
Suppose that $u\in W^{1,p}_{loc}(\R^N)\cap {\cal C} (\R^N)$  is a weak nonnegative  solution
of 
\begin{equation}\label{st2}
-\Delta_p u+h(x) u^q\leq 0,
\end{equation}
where $h$ satisfies (\ref{eqh}). If $q> p-1$, then $u\equiv 0$.
\end{thm}

\noindent {\bf Proof}: In a first step, we prove the result for 
$q< {N(p-1)+(\gamma+1)p\over N-p}$.

Let us multiply
(\ref{st2}) by
$u\zeta^\alpha$ where
$\zeta$ is some nonnegative cut-off function, supported in $B(0,2R)$ and equals
$1$ on $B(0,R)$, where $R$ is large enough to have $h(x) = a|x|^\gamma$ for $|x|\geq R$.  One may
choose such function with in addition
$|\nabla
\zeta|_\infty\leq {C\over R}$. After  integrating by parts, (\ref{st2}) becomes 
\begin{equation}\label{in1}
\int_{B_{2R}}|\nabla u|^p\zeta^\alpha+\alpha \int_{B_{2R}} |\nabla
u|^{p-2}
\nabla u.\nabla
\zeta u \zeta^{\alpha-1}+\int_{B_{2R}} h u^{q+1} \zeta^\alpha\leq 0.
\end{equation}
 Using Holders' inequality on the second integral one gets 
\begin{eqnarray}\label{in2}
|\int_{B(0,2R)\setminus B(0,R)} |\nabla u|^{p-2} \nabla u.\nabla
\zeta u \zeta^{\alpha-1}|&\leq& \alpha \left(\int_{B(0,2R)\setminus B(0,R)} |\nabla
u|^p\zeta^\alpha\right)^{p-1\over p}\cdot\nonumber \\
& &  \left(\int_{B(0,2R)\setminus B(0,R)} |\nabla \zeta|^p u^p
\zeta^{\alpha-p}\right)^{1\over p}.
\end{eqnarray}
Once more by Holders' inequality one has:
\begin{eqnarray}\label{in3}\int_{B(0,2R)\setminus B(0,R)} |\nabla \zeta|^p u^p
\zeta^{\alpha-p}&\leq& \left( \int_{B(0,2R)\setminus B(0,R)} u^{q+1} h
\zeta^\alpha\right)^{p\over q+1}\cdot\\
& & \left(\int_{B(0,2R)\setminus B(0,R)} |\nabla \zeta|^{p(q+1)\over q+1-p}
\zeta^{\alpha-p(q+1\over q+1-p)}h^{- p\over q+1-p}\right)^{1-{p\over q+1}}.\nonumber
\end{eqnarray}
Choosing $\alpha> {p(q+1)\over q+1-p}$, it is easy to see that
\begin{equation}\label{in4}
\int_{B(0,2R)\setminus B(0,R)} |\nabla \zeta|^{p(q+1)\over q+1-p}
\zeta^{\alpha-p(q+1\over q+1-p)}h^{- p\over q+1-p}\leq C R^{N-{\gamma p+p(q+1)\over q+1-p}}.
\end{equation}
Defining 
$$I_R = \int_{B_{2R}} |\nabla u|^p \zeta^\alpha $$
and 
$$J_R = \int_{B_{2R}} h \zeta^\alpha u^{q+1}.$$
Inserting (\ref{in2}), (\ref{in3}), (\ref{in4}) in (\ref{in1}), we have obtained that 
\begin{equation}\label{in5}
I_R+J_R\leq C\left(R^{N-{\gamma p+p(q+1)\over q+1-p}}\right)^{1-{p\over
q+1}} I_R^{p-1\over p} J_R^{1\over q+1}.
\end{equation}
Now choose $\beta=1+ {p\over (q+1)(p-1)}$ and 
$\delta:=\beta({p-1\over p})=\beta^\prime ({1\over q+1})$, then using Young
inequality, (\ref{in5}) becomes: 
\begin{eqnarray*}
I_R+J_R&\leq& C R^{\left(N-{\gamma p+p(q+1)\over q+1-p}\right){q+1-p\over
(q+1)p}} ({1\over \beta}I_R^{\delta}+{1\over \beta^\prime }J_R^{\delta})\\
&\leq &C R^{\left(N-{\gamma p+p(q+1)\over q+1-p}\right){q+1-p\over
(q+1)p}}(I_R+J_R)^{\delta}.
\end{eqnarray*}
 It is easy to see that $\delta <1$ when $p<q+1$,
and furthermore that 
$N-({\gamma p+p(q+1)\over q+1-p})< 0$ when 
$q<{(N+\gamma)(p-1)+p+\gamma\over N-p}$.

This concludes the proof, just let $R\rightarrow +\infty$ and then
$I_R+J_R\rightarrow 0$ which implies that $u\equiv 0$.

In a second step, we observe that if $u$ satisfies
(\ref{st2}), 
 then
$u^\alpha$ with $\alpha\geq 1$, is a solution of
$$-\Delta_p (u^\alpha)+ \alpha^{p-1} h
(u^\alpha)^{q+(\alpha-1)(p-1)\over \alpha }\leq 0.$$ More
precisely 
$$-\Delta_p (u^\alpha)+ \alpha^{p-1} h u^{q+(\alpha-1)(p-1)}\leq-\alpha^{p-1}
(\alpha-1)(p-1) u^{(\alpha-1)(p-1)-1}|\nabla u|^p.$$
This can be seen by taking $\tilde \varphi = \alpha^{p-1} u^{(\alpha-1)(p-1)}\varphi$,
$\varphi\in {\cal D} (\R^N)$, $\varphi\geq 0$ as test function in the equation 
$$-\Delta_p u+ hu^q\leq 0.$$
 As a consequence 
one has
$u\equiv 0$ as soon as
${q+(\alpha-1)(p-1)\over
\alpha }\leq {(N+\gamma)(p-1)+p+\gamma\over N-p}$. This will be always possible for
$\alpha$ large since 
$${\rm lim}_{\alpha\rightarrow +\infty} {q+(\alpha-1)(p-1)\over \alpha
} = p-1< {(N+\gamma)(p-1)+p+\gamma\over N-p}$$
 for $\gamma> -p$. Finally for any $q>p-1$ there
exists a power
$\alpha$ such that $u^\alpha\equiv  0$, hence $u\equiv 0$.

\bigskip
\noindent{\bf Acknowledgements}
This work was mainly done while the first author was visiting the Laboratoire d'Analyse,
G\'eom\'etrie et Mod\'elisation of the University of Cergy-Pontoise, she wishes to thank the
people of the laboratoire for the kind invitation and their welcome.

Isabeau Birindelli \\
Universit\`a di Roma ``La Sapienza'' \\
Piazzale Aldo moro, 5\\
00185 Roma, Italy \\
e mail: isabeau@mat.uniroma1.it

\medskip
\noindent Fran\c coise Demengel\\
Universit\'e de Cergy Pontoise,\\
Site de Saint-Martin, 2 Avenue Adolphe Chauvin\\
95302 Cergy Pontoise\\
e mail: Francoise.Demengel@math.u-cergy.fr
\end{document}